%
%
%
%
\input amstex
\documentstyle{amsppt}
\NoBlackBoxes

\topmatter
\title  
irreducible polynomials
in 
$F_q[x]$ versus
$1-qz$ in $\bold{C}[z]$\endtitle
\author Barry Brent\endauthor

\address 
4304 12th Avenue South, Minneapolis, MN 55407-3218 
\endaddress

\email barrybrent\@member.ams.org \endemail


\abstract 
Let $q$ be any prime power, let $F_q$ be the field with $q$ elements,
and let $N_n$ be the number of monic irreducible 
polynomials of degree $n$ in 
$F_q[x]$.  Let $z \in \bold{C}$ be close to zero.
We show that
$\prod_{n=1}^\infty (1-z^n)^{N_n} = 
1-qz$.
\endabstract


\endtopmatter

\document
\head 1. Introduction \endhead

In this note, we prove the following
\proclaim{Theorem}
For $q$ a prime power, let $F_q$ be the field with $q$ elements,
and let $N_n$ be the number of monic irreducible 
polynomials of degree $n$ in 
$F_q[x]$. Let $s>q$.
If $ z \in \bold{C}$ and $|z| < \frac 1s$, then

$$\prod_{n=1}^\infty (1-z^n)^{N_n} = 
1-qz.$$
\endproclaim

\head 2. Proof of the theorem \endhead
Divisor sums are 
taken over positive divisors throughout.  Let $\mu$ be the M\"obius function, and
let $N(a,n) = n^{-1}\sum_{d|n} \mu(n/d)a^d.$ It is clear that the 
number $N(a,n)$ is 
positive when $a>1$.
We have the following 
\proclaim{Lemma}  If $b > a > 1, z \in \bold{C}$ and $|z|<\frac 1b$, then
$$\prod_{n=1}^{\infty} (1 -z^n)^{N(a,n)} = 1 - az.$$
\endproclaim

\it 

\flushpar
Proof of the lemma \rm
First we show that the infinite product is absolutely convergent,
which is the same as showing that the series
$S = \sum_{n=1}^{\infty} N(a,n) | \log (1 - z^n)|$ converges. 
By  
M\"obius inversion, $$a^n = \sum_{d|n} d N(a,d), \tag{1}$$
so $N(a,n) \leq a^{n}/n$.
For small $|z|$ and large $n$ we have the existence of a positive constant $C$ such that
$|\log(1-z^n)| \le C|z^n|$ because 
$\lim_{u \rightarrow 0} \log (1-u)/u = 0$. Therefore the tail of 
the series $S$ is dominated for some $k$ by
the series $\sum_{n=k}^{\infty} C\frac {a^n}n | z^n|$, which, 
under our hypotheses, converges by the ratio test.
Hence the series $S$ converges
and the product
$\prod_{n=1}^{\infty} (1 -z^n)^{N(a,n)}$
converges absolutely.

Now clearly the series $\sum_{n=1}^{\infty} N(a,n) z^n$
also converges absolutely under our hypotheses.  We need these two
facts because we are going to apply the following 
result (due essentially to Euler):
an absolutely convergent infinite product $\prod_{n=1}^{\infty} 
(1-z^n)^{f(n)/n}$ with the property that
$\sum_{n=1}^{\infty} f(n)/n$ converges absolutely
expands as a power series $1 +\sum_{n=1}^{\infty}r(n)z^n$
with (for $n = 1, 2, 3, ... $) 
$$n  r(n) + \sum_{k=1}^n r(n-k)\sum_{\Sb  d|k \endSb }f(d) = 0  \tag{2}$$
\flushpar
(\cite{A}, Theorem 14.8).
 Specializing equation (2)
to the situation $f(n) = n N(a,n)$ and applying equation (1), we find that
$$ nr(n) + \sum_{k=1}^n r(n-k)a^k = 0.\tag{3}
$$
If we write $r(0) = 1$ and apply equation (3) twice with $n = 1$ and $n = 2$,
we find that $r(1) = -a$ and $r(2) = 0$. Then it is
easy to show that $r(n) = 0$ for all $n>1$ by induction. This finishes
the proof of the lemma.

If $q = p^n$, $p$ prime, then $N(q,n) =N_n$, 
the number of monic irreducible polynomials
in $F_q[x]$ 
(\cite{G},\cite{BeBl}, or, for the case 
$q = p$,
\cite{IR}).
So the theorem follows from the lemma.

\refstyle{C}

\enddocument